\documentclass[a4paper, 12pt]{article}
\usepackage{setspace} 
\usepackage{fullpage}
\usepackage{amsmath, amsthm, amssymb, mathrsfs, graphicx, subfigure}
\usepackage{ifthen, url}
\usepackage[usenames]{color}
\usepackage{hyperref}

\usepackage[sort&compress]{natbib}

\theoremstyle{plain}

\newtheorem{example}{Example}
\newtheorem{assumption}{Assumption}

\theoremstyle{remark}

\theoremstyle{definition}
\newtheorem{definition}{Definition}

\DeclareMathOperator*{\argmin}{arg\,min}
\makeatletter
\newcommand*{\defeq}{\mathrel{\rlap{%
                     \raisebox{0.3ex}{$\m@th\cdot$}}%
                     \raisebox{-0.3ex}{$\m@th\cdot$}}%
                     =}
\makeatother

\begin{document}

\title{Bayes Factor Consistency}

\author{Siddhartha Chib \\
John M. Olin School of Business, Washington University in St. Louis \\
\url{chib@wustl.edu} \\
\mbox{} \\
Todd A. Kuffner \\
Department of Mathematics, Washington University in St. Louis \\
\url{kuffner@math.wustl.edu}
}
\date{\today}

\maketitle

\onehalfspacing
\begin{abstract}
Good large sample performance is typically a minimum requirement of any model selection criterion. This article focuses on the consistency property of the Bayes factor, a commonly used model comparison tool, which has experienced a recent surge of attention in the literature. We thoroughly review existing results. As there exists such a wide variety of settings to be considered, e.g. parametric vs. nonparametric, nested vs. non-nested, etc., we adopt the view that a unified framework has didactic value. Using the basic marginal likelihood identity of Chib (1995), we study Bayes factor asymptotics by decomposing the natural logarithm of the ratio of marginal likelihoods into three components. These are, respectively, log ratios of likelihoods, prior densities, and posterior densities. This yields an interpretation of the log ratio of posteriors as a penalty term, and emphasizes that to understand Bayes factor consistency, the prior support conditions driving posterior consistency in each respective model under comparison should be contrasted in terms of the rates of posterior contraction they imply.
\medskip

\end{abstract}

\newpage
\doublespacing
\section{Introduction}

Bayes factors have long held a special place in the Bayesian inferential paradigm, being the criterion of choice in model comparison problems for such Bayesian stalwarts as Jeffreys, Good, Jaynes, and others. An excellent introduction to Bayes factors is given by \citet{KassRaftery:1995}. The meritorious reputation of the Bayes factor derives from its relatively good performance across key inference desiderata, including interpretability, Occam's razor, and an ability to choose the best model among those under comparison in large samples. One concrete definition of the latter property is consistency.

Informally, let $\mathcal{M}_{1}$ and $\mathcal{M}_{2}$ be the only two candidate statistical models, each specifying a set of distributions for the data and a prior distribution on this set. Assume \emph{a priori} that the models are assigned equal odds. The posterior under model $\mathcal{M}_{k}$, $k=1,2$, is given by
\begin{equation}
\label{intropost}
\{\text{posterior } | \mathcal{M}_{k}\}=\frac{\{\text{likelihood } | \mathcal{M}_{k}\} \times \{\text{prior } | \mathcal{M}_{k}\}}{\text{normalizing constant}}.
\end{equation}
The normalizing constant, which is simply the integral of $\{\text{likelihood } | \mathcal{M}_{k}\} \times \{\text{prior } | \mathcal{M}_{k}\}$ over the parameter space (which may be finite- or infinite-dimensional), is called the marginal likelihood under $\mathcal{M}_{k}$, denoted by $m(\text{data} | \mathcal{M}_{k})$.

The Bayes factor for comparing model $\mathcal{M}_{1}$ to model $\mathcal{M}_{2}$ is
\[
BF_{12}=\frac{m(\text{data} | \mathcal{M}_{1})}{m(\text{data} | \mathcal{M}_{2})}.
\]
A `large' value of $BF_{12}$ indicates support for $\mathcal{M}_{1}$ relative to $\mathcal{M}_{2}$, and a `small' value ($>0$) indicates support for $\mathcal{M}_{2}$ relative to $\mathcal{M}_{1}$. Bayes factor consistency refers to the stochastic convergence of $BF_{12}$, under the true probability distribution, such that $BF_{12} \rightarrow \infty$ if $\mathcal{M}_{1}$ is the best model, and $BF_{12} \rightarrow 0$ if $\mathcal{M}_{2}$ is the best model.

The marginal likelihood is generally not analytically tractable to compute, and hence must be approximated in practice. The conventional approaches to studying the large-sample properties of the Bayes factor involve studying the large-sample properties of the marginal likelihood, or a suitable approximation, under each model, or to derive suitable bounds for these quantities.

The literature has witnessed a surge in interest regarding Bayes factors for model comparison since the development of accurate asymptotic approximations in the Bayesian setting, and the advent of Markov Chain Monte Carlo (MCMC) methods for estimating the marginal likelihood. These breakthroughs enabled practitioners to overcome analytically intractable problems and brought the Bayes factor into everyday use. Key references include \citet{TierneyKadane:1986}, \citet{GelfandSmith:1990}, \citet{KassVaidyanathan:1992}, \citet{CarlinChib:1995}, \citet{Chib:1995}, \citet{Green:1995}, \citet{VerdinelliWasserman:1995}, \citet{MengWong:1996}, \citet{ChenShao:1997}, \citet{DKRW:1997}, \citet{ChibJeliazkov:2001}, \citet{HanCarlin:2001} and \citet{BasuChib:2003}.

In this article, we review existing consistency results through the lens of a simple decomposition. \citet{Chib:1995} developed a convenient approach for estimating the marginal likelihood using an identity found by a rearrangement of the posterior \eqref{intropost}. Namely,
\[
m(\text{data} | \mathcal{M}_{k})=\frac{\{\text{likelihood } | \mathcal{M}_{k}\} \times \{\text{prior } | \mathcal{M}_{k}\}}{\{\text{posterior } | \mathcal{M}_{k}\}}.
\]
We use this identity to write the natural logarithm of the Bayes factor as the sum of three log ratios, each of which can be studied in terms of their large sample behavior when evaluated at the same sequence of points. Specifically,
\[
\log BF_{12}=\log \frac{\{\text{likelihood } | \mathcal{M}_{1}\}}{\{\text{likelihood } | \mathcal{M}_{2}\}} + \log \frac{\{\text{prior } | \mathcal{M}_{1}\}}{\{\text{prior } | \mathcal{M}_{2}\}}-\log \frac{\{\text{posterior } | \mathcal{M}_{1}\}}{\{\text{posterior } | \mathcal{M}_{2}\}}.
\]
Whatever the target of inference may be in each model, such as a finite-dimensional parameter or a density, we argue that looking at the sample-size-dependent sequence of values of a consistent estimator for the target in each model, respectively, facilitates the application of a rich literature concerning likelihood and posterior asymptotics.

We first formally define the Bayes factor and its consistency property in the most general setting, and in the sequel make the concepts in each setting more precise. Throughout, we try to avoid technical details which, while nontrivial, are not essential to follow the arguments. We also focus only on the usual Bayes factor, rather than its variants. The practical aspects of the parametric and nonparametric settings are sufficiently different that they must be discussed separately, though we may accommodate both settings in our generic description of the problem below. For the technical details, the interested reader is directed to the relevant articles as indicated throughout. Recent monographs include \citet{GhoshRamamoorthi:2003}, \citet{HHMW:2010} and \citet{GineNickl:2015}.

\section{Problem Setting and Basic Argument}

\subsection{Notation and Definitions}
We borrow from notational conventions in the nonparametric literature, c.f. \citet{GLV:2008, Walker:2004a, HHMW:2010}. Let $\mathbf{y}^{(n)}\equiv \mathbf{y}$ be the observed values of a sequence of $n$ random variables $\mathbf{Y}^{(n)}=\{Y_{1}, \ldots, Y_{n}\}$ which are independent and identically distributed according to some distribution $P_{0}$ on a measurable space $(\mathcal{Y}, \mathcal{A})$, having density $p_{0}$ with respect to a suitable dominating measure $\mu$ on $(\mathcal{Y}, \mathcal{A})$. The joint distribution of $\mathbf{Y}^{(n)}$, which is the $n$-fold copy of $P_{0}$ denoted by $P_{0}^{n}$, is absolutely continuous with respect to a common measure $\mu^{n}$ on the sample space $\mathcal{Y}^{n}$. It is desired to choose the best model among two or more candidates, where `best' intuitively means that it most closely resembles the unknown true distribution of the data generating process, $P_{0}^{n}$. The `best' model may be within some class of parametric models, semiparametric models, or fully nonparametric models. The simplest case is when the `best' model corresponds to the true model, i.e. when the truth is contained in one of the models under consideration.

For sample size $n$, a generic model is denoted by $M_{n,\alpha}=\{\mathcal{F}_{n,\alpha}, \Pi_{n,\alpha}, \lambda_{n, \alpha}\}$, where $\mathcal{F}$ is a set of $\mu$-probability densities on $(\mathcal{Y}, \mathcal{A})$ equipped with a $\sigma$-algebra ensuring that the maps $(y, p) \mapsto p(y)$ are measurable. The models are indexed by $\alpha \in A_{n}$, where for each $n \in \mathbb{N}$, the index set $A_{n}$ is countable. The prior distribution $\Pi_{n,\alpha}$ is a probability measure on $\mathcal{F}_{n,\alpha}$, and $\lambda_{n, \alpha}$ is a probability measure on $A_{n}$. The assumption that the index set is countable for every $n$ essentially removes some technical problems which would arise if all models had zero prior probability, which would occur if there were infinitely many possible models.

We allow that the sets $\mathcal{F}_{k}$ and/or $\mathcal{F}_{l}$ may be of any general form. That is, $M_{k}$ and $M_{l}$ could be any combination of parametric, semiparametric or nonparametric models. For example, in a parametric model $M$, $\mathcal{F}=\{p_{\boldsymbol{\theta}}^{n}(\mathbf{y}), \boldsymbol{\theta} \in \boldsymbol{\Omega} \subset \mathbb{R}^{d}, \pi(\boldsymbol{\theta})\}$, where $p_{\boldsymbol{\theta}}^{n}(\mathbf{y})$ denotes the density of $\mathbf{y}$ with respect to $\mu^{n}$ under model $M$ prescribing the set $\mathcal{F}$.

For a generic model $M$, the overall prior is a probability measure on the set of probability densities,
\begin{equation}
\label{prior}
\Pi_{n}=\sum_{\alpha \in A_{n}} \lambda_{n,\alpha}\Pi_{n,\alpha}.
\end{equation}
The posterior distribution of a model index $B \subset A_{n}$, given this prior distribution, is the random measure
\begin{align*}
\label{posterior}
\Pi_{n}(B | y_{1}, \ldots, y_{n})&=\frac{\int_{B}\prod_{i=1}^{n}p(y_{i})d\Pi_{n}(p)}{\int \prod_{i=1}^{n}p(y_{i})d\Pi_{n}(p)} \\
&= \frac{\sum_{\alpha \in A_{n}}\lambda_{n,\alpha}\int_{p \in \mathcal{P}_{n,\alpha}: p \in B}\prod_{i=1}^{n}p(y_{i})d\Pi_{n,\alpha}(p)}{\sum_{\alpha \in A_{n}} \lambda_{n,\alpha} \int_{p \in \mathcal{P}_{n,\alpha}}\prod_{i=1}^{n}p(y_{i})d\Pi_{n,\alpha}(p)},
\end{align*}
and the marginal likelihood is defined as
\begin{equation}
\label{marglik}
m(\mathbf{y} | M)=\frac{\sum_{\alpha \in A_{n}}\lambda_{n,\alpha}\int_{p \in \mathcal{P}_{n,\alpha}: p \in B}\prod_{i=1}^{n}p(y_{i})d\Pi_{n,\alpha}(p)}{\Pi_{n}(B | y_{1}, \ldots, y_{n})}.
\end{equation}
For clarity, note that the given, fixed true density $p_{0}$ corresponds to a `best' index element, say $\beta_{n} \in A_{n}$, in the sense that $\beta_{n}$ refers to the model which is closest to the true model according to a chosen measure. In the above, the set $B \subset A_{n}$ of indices could simply be the single element $\beta_{n}$, or a collection of models.

Consider two candidate models $M_{k}$ and $M_{l}$, prescribing sets of density functions $\mathcal{F}_{k}$ and $\mathcal{F}_{l}$, with associated prior distributions on these sets of density functions $\Pi_{k}$ and $\Pi_{l}$, and let $m(\mathbf{y} | M_{k})$ and $m(\mathbf{y} | M_{l})$ denote the respective marginal likelihoods.  The Bayes factor is defined by
\begin{equation}
\label{eq:BFdef}
BF_{kl}=\frac{m(\mathbf{y} | M_{k})}{m(\mathbf{y} | M_{l})}  = \frac{\lambda_{n,k}\int \prod_{i=1}^{n}p(y_{i})\Pi_{n,k}(p)}{\lambda_{n,l}\int \prod_{i=1}^{n}p(y_{i})\Pi_{n,l}(p)} .
\end{equation}

The prior probabilities on the models, $\lambda_{n, \alpha}$, do not affect the consistency arguments, and thus there is no loss of generality in assuming $\lambda_{n, k}=\lambda_{n, l}$, so that these quantities may be ignored in what follows.

Consistency of the Bayes factor refers to stochastic convergence of the quantity (\ref{eq:BFdef}).
\begin{definition}[Conventional Bayes Factor Consistency]
\label{def:BF}
The Bayes factor for comparing $M_{k}$ and $M_{l}$, $BF_{kl}=m(\mathbf{y} | M_{k})/m(\mathbf{y} | M_{l})$, is consistent if:
\begin{enumerate}
\item[(i)] $BF_{kl} \rightarrow_{p} 0$ (or $\log BF_{kl} \rightarrow_{p} -\infty$) when $M_{l}$ contains the true model ($p_{0}^{n} \in \mathcal{F}_{n,l}$); and
\item[(ii)] $BF_{kl} \rightarrow_{p} \infty$ (or $\log BF_{kl} \rightarrow_{p} \infty$) when $M_{k}$ contains the true model ($p_{0}^{n} \in \mathcal{F}_{n,k}$).
\end{enumerate}
The probability measure associated with these convergence results is the one associated with the infinite-dimensional product measure corresponding to the true distribution, the $n$-fold product measure $P_{0}^{n}$ as $n \rightarrow \infty$. When both relations hold with probability one, the Bayes factor is said to be almost surely consistent. In all stochastic convergence statements in this paper, we are considering $n \rightarrow \infty$. We further note that such convergence statements are pointwise and not uniform.
\end{definition}
It will often be convenient to work with the natural logarithm of the Bayes factor, referred to by I.~J. Good as the \emph{weight of evidence}. As noted in the definition, and pointed out by \citet{DBG:2012} and \citet{ChakrabartiGhosh:2011}, we must be careful about which probability measure is associated with the stochastic convergence statement. Moreover, consistency of the Bayes factor in the sense of Definition~\ref{def:BF} requires that one of the models being considered contains the true model. Definition~\ref{def:NP_BFdef} accommodates the more general setting.

Bayes factor consistency is not the same as model selection consistency. The reason this is often referred to as model selection consistency is that this ensures the sequence of posterior probabilities of the true model will converge to one, at least in fixed-dimensional settings; see \citet{LPMCB:2008}, \citet{CGMM:2009} and \citet{ShangClayton:2011}. This is what is more conventionally thought of as model selection consistency \citep{FLS:2001}. It has been emphasized \citep{MGC:2015} that the Bayesian and frequentist notions of model selection consistency do not necessarily agree; c.f. \citet{Shao:1997} where model selection consistency means convergence in probability of the selected model to the submodel which minimizes the mean squared prediction error.

Note that if a prior is improper, then the marginal likelihood is also improper. In that case, the ratio (\ref{eq:BFdef}) cannot be interpreted. We consider here only proper priors for which the marginal likelihoods and Bayes factors are well-defined. There is a growing literature concerning alternative measures of evidence which allow for improper priors \citep{BergerPericchi:2001}, e.g. posterior Bayes factors \citep{Aitkin:1991}, fractional Bayes factors \citep{OHagan:1995} and intrinsic Bayes factors \citep{BergerPericchi:1996}. Other notable proposals for measuring evidence, as alternatives to the usual Bayes factor considered here, include posterior likelihood ratios \citep{Dempster:1973, SmithFerrari:2014}, test martingales \citep{SSVV:2011} and relative belief ratios \citep{Evans:2015}. A recent overview of objective Bayesian approaches can be found in \citet{BBFG:2012}.

\subsection{A Unified Framework for Analysis}

In this general setup, following the basic marginal likelihood identity (BMI) of \citet{Chib:1995}, we have that
\begin{equation}
\label{ChibBMI}
m(\mathbf{y} | M_{n, \alpha})= \frac{ p^{n}(\mathbf{y} | \mathcal{F}_{n, \alpha})\Pi_{n,\alpha}(\mathcal{F}_{n, \alpha})}{\Pi_{n, \alpha}(\mathcal{F}_{n, \alpha} | \mathbf{y})} .
\end{equation}
\noindent
Using this identity, the natural logarithm of the Bayes factor for comparing $M_{k}$ and $M_{l}$ may be expressed as
\begin{equation}
\label{logBF}
\log BF_{kl}=\log \frac{p^{n}(\mathbf{y} | \mathcal{F}_{n, k})}{p^{n}(\mathbf{y} | \mathcal{F}_{n, l})} + \log \frac{\Pi_{n, k}(\mathcal{F}_{n, k})}{\Pi_{n,l}(\mathcal{F}_{n,l})} -\log \frac{\Pi_{n,k}(\mathcal{F}_{n, k} | \mathbf{y})}{\Pi_{n, l}(\mathcal{F}_{n,l} | \mathbf{y})}.
\end{equation}
To establish stochastic convergence as $n \rightarrow \infty$ under probability law $P_{0}^{n}$, when $p_{0}^{n} \in \mathcal{F}_{n, k}$ or $p_{0}^{n} \in \mathcal{F}_{n,l}$, we examine each of the three terms appearing on the r.h.s. of \eqref{logBF} separately. We emphasize that no single technique will work across all model comparison problems; the tools needed for analysis of each term in this decomposition will depend on the nature of the model comparison problem.

The first term is a type of log likelihood ratio. It is convenient for this term to be bounded in $P_{0}^{n}$-probability and, hence, to be $O_{p}(1)$. This happens, for example, if this quantity converges in distribution. The arsenal of tools for this term includes everything from classical likelihood theory to generalized likelihood ratio tests \citep{FZZ:2001} and other variants for the semiparametric setting. However, we emphasize that the goal is not simply to make assumptions which ensure convergence in distribution, but to also understand settings for which this would fail to happen, as such failure will impact Bayes factor consistency.

The second term is desired to be $O(1)$, as we want the prior in each model to have a negligible effect on the posterior for large samples. If the priors are continuous and bounded, as they will be when proper, this term will be bounded. This term can be problematic in some cases, however, including when the dimension of the parameter space grows with the sample size.

The final term involves a log ratio of posteriors. Consistency requires that this converge in $P_{0}^{n}$-probability to $-\infty$ if $p_{0} \in \mathcal{F}_{n,l}$ and to $\rightarrow \infty$ in $P_{0}^{n}$-probability if $p_{0}^{n} \in \mathcal{F}_{n,k}$. Essentially this just means that the correct model has a faster rate of posterior contraction, as a function of the available sample size, when compared to an incorrect or more complex alternative. This can often be established by existing theorems, provided one is willing to make suitable assumptions on the concentration of the priors, and the size of the model space. We review such conditions below.

An attractive feature of this framework is that the third term acts like a penalty term, and it is this term which typically drives the consistency result. This observation leads to interesting questions about when the penalty is sufficiently large to prevent the selection of an incorrect model, and to prevent overfitting. Intuitively, if two models are close in some sense but only one of them is correct, consistency requires some condition ensuring that the posterior in the correct model contracts faster, and that it is possible to distinguish the models asymptotically. The latter condition can be roughly thought of as saying the frequentist notion of the `null' and `alternative' hypotheses are well-separated in an appropriate sense. This is related to the existence of uniformly exponentially consistent tests. Moreover, if the models are nested, and both contain the truth, the penalty term must be sufficiently large to ensure the smaller model is chosen. In the classical nested parametric setting, the larger model, by virtue of having a higher-dimensional parameter space, will have prior mass that is more spread out. Thus if the true density is contained in a smaller model, with a more concentrated prior mass around the true density, the posterior in the smaller model will contract faster than the more complex model. In more general, non-nested settings with nonparametric models, a similar intuition will hold.

If the true $P_{0}$ is not contained in the set of candidate models, it will certainly not be possible for any model selection procedure to select the correct model. Many authors argue that it is unrealistic to assume the true distribution of the data generating process is exactly specified by one of the candidate models, and instead use the qualifier `pseudo-true' to mean the `best' approximating model among those being considered. While we sympathize with such authors and use this convention when appropriate, we often make no distinction in what follows between true and pseudo-true models (and parameter values). This issue is most relevant when both models are misspecified and not well-separated, so that it is difficult to construct a consistent sequence of tests. Moreover, while consistency in the usual sense is not possible when none of the candidate models contain the true distribution, it will still be possible to define some notion of consistency in the sense of choosing the model which is closest to the true distribution in a relative entropy sense; c.f. Definition~\ref{def:NP_BFdef}.

\subsection{Advantages of This Approach}

Since \eqref{ChibBMI} is an identity, then it holds for any density in $\mathcal{F}$, the set of densities being considered. It will also hold for every sequence of estimates of the density (or parameter, in the parametric setting). Thus, to establish consistency of the Bayes factor, it is sufficient to show convergence in $P_{0}^{n}$-probability of the Bayes factor for any sequence of estimates which are the values of consistent estimators for the unknown densities or parameters in the respective models being compared.

A second advantage is that we examine the asymptotic behavior of each component, rather than the aggregate asymptotic behavior of the marginal likelihood, and in this sense we learn more about the asymptotic behavior of important quantities. We are able to exploit a richer literature to study likelihood ratios and posterior contraction rates separately, and this can help identify important unanswered questions or problems for which there is scope to improve existing results. Compared to existing proofs, which rely on results which relate the prior support conditions to bounds on the marginal likelihood, we are directly considering how the prior support conditions affect the behavior of the posterior in each model, which is conceptually appealing.

\subsection{Model Comparison Types}

When comparing two models, it is possible that both, only one or neither of the models contains the true distribution. When a model contains the true distribution, it is said to be correctly specified; otherwise it is said to be misspecified. We can further classify the types of model comparison as: (a) both models are parametric (\S~\ref{sec:parpar}); (b) one model is parametric, the other is nonparametric (\S~\ref{sec:parnonpar}); (c) one model is parametric, the other is semiparametric (\S~\ref{sec:parsemi}); (d) both models are nonparametric (\S~\ref{sec:discussion}); (e) both models are semiparametric (\S~\ref{sec:discussion}); (f) one model is nonparametric, the other is semiparametric (\S~\ref{sec:discussion}).

Thus there are six frameworks of model comparison and, within any of these, one could consider misspecified models, non i.i.d. observations, or other `non-regular' settings.  Much of the Bayes factor literature makes distinctions between the situations where the models are nested, overlapping or non-nested. These classifications may be roughly understood in the parametric setting as follows. Model 1 is nested in Model 2 if it is possible to obtain Model 1 by some restriction on the parameters of Model 2. This is the case most commonly studied in classical frequentist hypothesis testing for linear models; the relationship between the parameter spaces and likelihoods facilitates the derivation of the asymptotic distribution of the likelihood ratio statistic under the null model. Non-nested models can be either strictly non-nested or overlapping. Overlapping models are those for which neither model is nested in the other, but that there still exists some choice of parameter values for each model such that they yield the same joint distribution for the data. Strictly non-nested models do not allow for this latter possibility. In the nonparametric setting, borrowing from the literature on adaptation, the distinction is in terms of coarser and smoother models. More recently, \citet{GGV:2000} and \citet{GLV:2008} have adopted the terminology of bigger and smaller models, with a precise definition of model complexity given in \citet{GGV:2000}. We also note that, while the classical large sample evaluation of Bayes factors considers the dimension of the model to be fixed, there is a growing literature concerning the asymptotic behavior of Bayes factors as the model dimension grows with the sample size. We reference this literature throughout.

We focus on the `regular' versions of each of the above frameworks, with comments and references in appropriate places regarding extensions to common alternative settings. While we have tried to incorporate as many of the recent developments as possible, the size of the Bayes factors literature renders any attempt at being comprehensive as an exercise in futility, and our omission of certain elements should not be interpreted as a judgment that these contributions are less important.

We caution the reader that in any particular setting, there will be model-specific assumptions required to ensure that the model comparison problem is well-defined. For example, in the regression settings we mention below, there would be some additional assumptions on the matrix of predictor variables, the error distribution and the joint distribution of the predictors and errors to ensure the consistency of a suitable estimator for the unknown mean regression function or its parameters, and hence also for a consistent estimator of the density of interest. These assumptions are crucial, but we do not dwell on them here. Instead, we focus only on the aspects of the problem which are unique to the study of Bayes factor asymptotics.

\section{Frameworks, Concepts and Connections}
We review some key definitions and concepts needed for large-sample model comparison across the parametric, semiparametric and nonparametric frameworks. Posterior consistency and rates of contraction are also considered. A running example of mean regression modeling is introduced.

\subsection{Parametric Framework}

Let $\mathbf{Y}^{(n)}=\{Y_{1}, \ldots, Y_{n}\}$ be a sequence of random variables from a population with some density indexed by parameter vector $\boldsymbol{\theta} \in \Omega$; the $Y_{i}$ are assumed to be conditionally independent and identically distributed, given $\boldsymbol{\theta}$.  The parameter sets considered are subsets of $\mathbb{R}^{d}$, where $d$ is finite and fixed, in particular, $d<n$ whenever estimation is considered for fixed $n$. A model, $M_{k}$, $k=1, \ldots, K$ (where $K \in \mathbb{N}$, hence a countably infinite set) consists of a parameter space $\Omega_{k}$, a density $f_{k}(\mathbf{y} | \boldsymbol{\theta}_{k} \in \Omega_{k}, M_{k})$ and a prior density $\pi_{k}(\boldsymbol{\theta}_{k} | M_{k})$ for $\boldsymbol{\theta}_{k} \in \Omega_{k}$.  The vector of observed values $\mathbf{y}^{(n)}=(y_{1}, \ldots, y_{n})$ are the realizations of $\mathbf{Y}^{(n)}$, though we suppress the superscript on $\mathbf{y}^{(n)} \equiv \mathbf{y}$ for notational convenience. The model space we will consider is $\mathcal{M}=\cup_{k=1}^{K}M_{k}$, a countable union of models.

A Bayesian constructs a prior distribution $\Pi$, which expresses beliefs about the parameter. Combining the prior with the observations yields the posterior distribution
\begin{equation}
\label{posterior_parametric}
\Pi_{n}(\boldsymbol{\theta} |\mathbf{y})=\frac{\pi(\boldsymbol{\theta})f(\mathbf{y} | \boldsymbol{\theta})}{\int_{\boldsymbol{\Omega}}\pi(\boldsymbol{\theta})f(\mathbf{y} | \boldsymbol{\theta})d\boldsymbol{\theta}}=\frac{\pi(\boldsymbol{\theta})f(\mathbf{y} | \boldsymbol{\theta})}{m(\mathbf{y})}
 \end{equation}
with the prior $\pi(\boldsymbol{\theta})$ representing the density corresponding to the prior probability measure $\Pi$. The marginal likelihood under model $M_{k}$ is
\begin{equation}
\label{ML_parametric}
m(\mathbf{y} | M_{k})=\int f(\mathbf{y} | \boldsymbol{\theta}_{k}, M_{k})\pi(\boldsymbol{\theta}_{k} | M_{k})d\boldsymbol{\theta}_{k}.
\end{equation}

\begin{example}[Linear Regression]
Let $\{Y_{i}\}_{i=1}^{n}=\{(Z_{i}, X_{ij})\}_{i=1}^{n}$, $j=1, \ldots, p$, so that $\mathbf{Y}$ is a matrix of response-covariate pairs, where $\mathbf{Z}=\{Z_{1}, \ldots, Z_{n}\}$ is an $n$-dimensional vector of response variables and $\mathbf{X} \in \mathbb{R}^{n \times p}$ is the design matrix. It is assumed that
\begin{equation}
Z_{i} = X_{ij}\boldsymbol{\beta}+\varepsilon_{i}
\end{equation}

where $\boldsymbol{\beta} \in \mathbb{R}^{p}$ is an unknown parameter vector. We typically assume that the $\epsilon_{i}$ are i.i.d. $\mathcal{N}(0, \sigma^{2})$ with $\sigma >0$, and thus $\mathbf{Z}| \mathbf{X} \sim \mathcal{N}(\mathbf{X}\boldsymbol{\beta}, \sigma^{2}I_{n})$. It is typically assumed that an intercept term is included, so that the first column of $\mathbf{X}$ is a vector of $1$s. It is also conventional for illustrative purposes to assume that the $X_{ij}$ are fixed. A Bayesian model specifies a prior for $\boldsymbol{\beta} \sim \Pi_{1}$, independent of $\sigma^{2}$, and could either treat $\sigma^{2}$ as known or specify a prior $\sigma^{2} \sim \Pi_{2}$.
\end{example}

Complete specification of the Bayesian model comparison procedure would also require a prior distribution over the space of models; this would be necessary to calculate the posterior odds from the Bayes factor. This will not be necessary for our purposes. The priors for the parameter are assumed to be proper; use of improper priors in conjunction with Bayes factors for model selection can be problematic (e.g. Lindley's paradox). We refer readers to \citet{Robert:1993} and \citet{VillaWalker:2015} for discussion of Lindley's paradox. \citet{GPS:2005} point out that the paradox disappears if the Bayesian and frequentist asymptotic frameworks are the same, i.e. both using Bahadur asymptotics (fixed alternatives) or Pitman (contiguous) alternatives. When improper priors are used, such as those arising from objective Bayesian analysis, there have been several proposals to alter the standard model comparison framework. \citet{Dawid:2011}, suggested that using the posterior odds instead of the Bayes factor can overcome the problem, while \citet{OHagan:1995}, \citet{DawidMusio:2015} and others have proposed variants such as fractional or intrinsic Bayes factors.

\subsection{Nonparametric Framework}

Let $\mathbf{Y}^{(n)}=\{Y_{1}, \ldots, Y_{n}\}$ be a sequence of random variables, which are assumed to be independent and identically distributed according to a true distribution $F_{0}$, having true density $f_{0}$ with respect to Lebesgue measure. In the parametric setting, estimation is concerned with the parameter $\boldsymbol{\theta}$ of an (assumed) known distribution, while in the nonparametric setting, estimation is concerned with a density function. A Bayesian constructs a prior distribution $\Pi$ on the set of densities $\mathcal{F}$ prescribed by the model, which expresses beliefs about the location of the true density. The posterior distribution of a set of densities $A$ in the set $\Omega$ of all densities with respect to Lebesgue measure is given by
 \begin{equation}
 \label{posterior_nonparametric}
 \Pi_{n}(A)=\frac{\int_{A} f(\mathbf{y})\Pi(df)}{\int_{\Omega} f(\mathbf{y})\Pi(df)}= \frac{\int_{A} f(\mathbf{y})\Pi(df)}{m(\mathbf{y})},
 \end{equation}
 where $m(\mathbf{y})$ again denotes the marginal likelihood.

Given two densities $f$ and $g$, which are both absolutely continuous with respect to the same dominating measure $\mu$ over a set $\mathcal{S}$, define the Kullback-Leibler divergence of $g$ from $f$ as $d_{KL}(f,g) \equiv d_{KL}(g \| f)=\int_{\mathcal{S}} g \log (g/f)d\mu$. Conventionally the `best' density would be $g$, so that one is speaking of the divergence of some proposed distribution $f$ from the `best' density. The simplest case is when $g \equiv p_{0}$. A Kullback-Leibler neighborhood of the density $g$, of size $\epsilon >0$, is defined by
\begin{equation}
\label{KL_neighbor}
N_{g}(\epsilon)=\left\{ f : \int g(x) \log\frac{g(x)}{f(x)}dx < \epsilon \right\}.
\end{equation}

\begin{definition}[Kullback-Leibler Property]
A prior distribution $\Pi$ over the space of densities is said to possess the Kullback-Leibler property if
\begin{equation}
\label{KL_property}
\Pi \{ f :  d_{KL}(f,g) < \epsilon \} > 0
\end{equation}
for all $\epsilon > 0$ and for all $g$.
\end{definition}
Expositions of the Kullback-Leibler property, with examples, are found in \citet{WuGhosal:2008}, \citet{WDL:2004}, \citet{PetroneWasserman:2002}, \citet{GGR:1999a}, and \citet{BSW:1999}.

Each model $M_{\alpha}$ under consideration contains an element $f_{\alpha}$ which is such that
\[
f_{\alpha} \defeq \argmin _{f \in \mathcal{F}_{\alpha}} d_{KL}(f, p_{0}).
\]
The element $f_{\alpha}$ is the closest member of model $M_{\alpha}$ to the true density $p_{0}$. If the model contains the truth, then $d_{KL}(f_{\alpha}, p_{0})=0$. Otherwise the `best' density in $M_{\alpha}$ is simply defined as the closest density to the true density in the Kullback-Leibler sense. Consider two models $M_{k}$ and $M_{l}$. Each model contains a density, respectively $f_{k}$ and $f_{l}$, defined in the sense above.
 \begin{definition}[Pragmatic Bayes Factor Consistency \citep{WDL:2004}]
 \label{def:NP_BFdef}
A pragmatic version of Bayes factor consistency is said to hold if
\begin{enumerate}
\item[(i)] $BF_{kl} \rightarrow 0$ almost surely (or $\log BF_{kl} \rightarrow -\infty$ a.s.) when $d_{KL}(f_{l}, p_{0}) < d_{KL}(f_{\alpha}, p_{0})$ for all $\alpha \neq l$; and
\item[(ii)] $BF_{kl} \rightarrow \infty$ almost surely (or $\log BF_{kl} \rightarrow \infty$ a.s.) when $d_{KL}(f_{k}, p_{0}) < d_{KL}(f_{\alpha}, p_{0})$ for all $\alpha \neq k$.
\end{enumerate}
 \end{definition}
Again, the probability measure associated with the stochastic convergence statement is the infinite product measure $P_{0}^{\infty}$. This says that the Bayes factor will asymptotically choose the model containing the density which is closest, in the Kullback-Leibler sense, to the true model, among all possible densities contained in the models under comparison.

A fundamentally important result related to Definition~\ref{def:NP_BFdef} is given in Theorem 1 of \citet{WDL:2004}, which concerns the general setting where it is possible that neither model contains the true density. Formally, Bayesian models are characterized by a prior $\Pi$, and associated with a value $\delta_{\Pi} \geq 0$, which is such that $\Pi \{f : d_{KL}(f, p_{0}) < d \} > 0$ only and for all $d > \delta_{\Pi}$. Let $\delta_{\alpha} \geq 0$ be the value associated with $\Pi_{\alpha}$. \citet{WDL:2004} show that $BF_{KL} \rightarrow \infty$ almost surely if and only if $\delta_{k} < \delta_{l}$. Moreover, if one model has the Kullback-Leibler property while the other does not, then the Bayes factor will asymptotically prefer the model possessing the Kullback-Leibler property. This was also shown in \citet{DassLee:2004} for the setting of a testing a point null hypothesis against a non-parametric alternative with a prior possessing the Kullback-Leibler property. The proof in the latter paper, however, does not generalize to general null models; c.f. \citet{GLV:2008}. A related result due to \citet{CSS:2016} is that the less misspecified model, in the Kullback-Leibler sense, will asymptotically yield a larger marginal likelihood with probability tending to one. This is a result about model selection consistency when using the relative (pairwise) marginal likelihood as a selection criterion. 

When two or more models under consideration have the Kullback-Leibler property, stronger conditions are needed to ensure Bayes factor consistency. These are explicitly discussed in \citet{GLV:2008} and \citet{MRM:2009}. We discuss these conditions in \S~\ref{sec:parnonpar}.

 \begin{example}[Nonparametric Mean Regression]
 \label{ex:NPreg}
 A nonparametric regression model specifies
 \begin{equation}
 \label{NPreg}
 Z_{i}= r(X_{ij})+ \varepsilon_{i},
 \end{equation}
where $r(\cdot)$ is an unknown function, and the $\varepsilon_{i}$ are i.i.d. with mean zero. Again, for illustrative purposes it is typically assumed that the $X_{ij}$ are fixed.
 \end{example}

 \subsection{Semiparametric Framework}

In the semiparametric framework, it is assumed that a model has both parametric and nonparametric components. Then the semiparametric prior has two ingredients, $\Pi_{\text{Par}}$ for the parametric part and $\Pi_{\text{NP}}$ for the nonparametric part. The overall prior is given by $\Pi_{\text{Par}} \times \Pi_{\text{NP}}$.

Within the mean regression modeling setting, there are two common manifestations of semiparametric models. The first is a partially linear model, and the second is found in \citet{KunduDunson:2014}.
\begin{example}[Partially Linear Model with Known Error Distribution]
\begin{equation}
\label{PartiallyLinear}
Z_{i}= X_{ij} \boldsymbol{\beta} + r(X_{ij})+ \varepsilon_{i},
\end{equation}
where $r(\cdot)$ is an unknown function in an infinite-dimensional parameter space, and $X_{ij}\boldsymbol{\beta}$ is a linear component.
\end{example}

A semiparametric prior is $\Pi = \Pi_{\boldsymbol{\beta}} \times \Pi_{r}$. See Example~\ref{partlinear_compare} in \S~\ref{sec:parsemi} for more details. A simple version of the above example would be to modify Example~\ref{ex:NPreg}.  Assume $\varepsilon_{i} \sim \mathcal{N}(0, \sigma^{2})$ given $\sigma$, that $r(\cdot) \sim \Pi_{1}$ independent of $\sigma$ and $\boldsymbol{\varepsilon}$, and that $\sigma \sim \Pi_{2}$ where $\Pi_{1}$ and $\Pi_{2}$ are prior distributions. Then, formally, this would be a semiparametric model rather than fully nonparametric.

\begin{example}[Linear Model with Unknown Error Distribution]
\label{Lin_UnknownErr}
\begin{equation}
Z_{i}=X_{ij}\boldsymbol{\beta}+\varepsilon_{i}, \;\;\; \varepsilon_{i} \sim q(\cdot)
\end{equation}
where $q(\cdot)$ is an unknown residual density, about which parametric assumptions are avoided.
\end{example}

\subsection{Posterior Consistency}

The random measure $\Pi_{n}(\cdot | \mathbf{y})$, which is the posterior distribution, is said to be consistent for some fixed measure $P_{0}$ if it concentrates in arbitrarily small neighborhoods of $P_{0}$, either with probability tending to 1 or almost surely, as $n \rightarrow \infty$. Posterior consistency implies consistency of Bayesian estimates in the frequentist sense. See \citet{Barron:1986} for more on the definitions of weak, strong and intermediate forms of consistency. Parametric Bayesian asymptotics are reviewed in \citet[Ch. 1]{GhoshRamamoorthi:2003}. More recent contributions in Bayesian nonparametric asymptotics, where one may find other important references, include \citet{Castillo:2014} and \citet{GineNickl:2015}. We review here some essential results.

\citet{Schwartz:1965} established weak consistency of the posterior distribution under the condition that the prior has the Kullback-Leibler property.
\begin{definition}[Weak Consistency]
\label{def:WeakConsistency}
Let $F_{0}$ and $F$ be the cumulative distribution functions for the true density $p_{0}$ and a generic density $f$. Take any metric $w$ on the cumulative distribution functions $F_{0}$ and $F$ for which convergence in $w$ is equivalent to convergence in distribution. Define a weak neighborhood of $p_{0}$ to be $W=\{f : w(F_{0}, F) < \delta\}$ for $\delta>0$. A posterior distribution $\Pi_{n}(\cdot | \mathbf{y})$ is said to be weakly consistent if for almost all sequences under $p_{0}$, $\Pi_{n}(W | \mathbf{y}) \rightarrow 1$ for all weak neighborhoods of $p_{0}$.
\end{definition}
After \citet{Doob:1949} pioneered the notion of posterior consistency in the weak topology on the space of densities, as in the definition above, \citet{Schwartz:1965} established the Kullback-Leibler property as an important criterion for demonstrating weak consistency. However, \citet{DiaconisFreedman:1986a, DiaconisFreedman:1986b} have shown that priors satisfying this property in weak neighborhoods may not yield weakly consistent posteriors. \citet{Freedman:1963} and \citet{KimLee:2001} have also given examples of inconsistency. After the Diaconis-Freedman critique and insightful analysis by Andrew Barron, among others, the focus in the literature then shifted to establishing sufficient conditions for strong consistency, specifically Hellinger consistency; see \citet{Wasserman:1998}, \citet{GGV:2000}, \citet{ShenWasserman:2001}, \citet{WalkerHjort:2001}, \citet{Walker:2003} and \citet{Walker:2004a, Walker:2004b}. Another commonly used distance is the $L_{1}$ distance. By Hellinger consistency, we mean convergence according to the Hellinger metric on the set of densities. This metric induces the Hellinger topology on the set of densities.

\begin{definition}[Hellinger Consistency]
\label{def: HellingerConsistency}
For any pair of probability measures $P$ and $Q$ on a measurable space $(\mathcal{Y}, \mathcal{A})$ with corresponding densities $p$ and $q$ with respect to a dominating measure $\mu$, the Hellinger distance is defined as
\begin{equation}
\label{hellinger}
d_{H}(p, q)=\left\{\int \left(\sqrt{p}-\sqrt{q}\right)^{2}d\mu \right\}^{1/2}.
\end{equation}
Furthermore, define a Hellinger neighborhood of the true density to be $S_{\delta}=\{f : d_{H}(p_{0}, f) < \delta \}$. A sequence of posteriors $\{\Pi_{n}(\cdot | \mathbf{y})\}$ is Hellinger consistent for $p_{0}$, often called strongly consistent, if $\Pi_{n}(S_{\delta} | \mathbf{y}) \rightarrow 1$ almost surely in $P_{0}^{\infty}$ probability as $n \rightarrow \infty$.
\end{definition}

Posterior consistency in nonparametric and semiparametric problems is still an active research area, though some important questions have been answered. Some relatively recent and important contributions include \citet{BSW:1999}, \citet{Walker:2004b}, \citet{ShenWasserman:2001}. Particularly lucid overviews of nonparametric asymptotics are given by \citet{Ghosal:2010}, \citet{MartinHong:2012} and \citet[\S 7.4]{GineNickl:2015}. Pseudo-posterior consistency was considered in \citet{WalkerHjort:2001}.

In the semiparametric setting, letting $\boldsymbol{\theta}$ and $q$ respectively denote the finite-dimensional parametric and inifinite-dimensional nonparametric components, posterior consistency can refer to (i) consistency at the pair $(\boldsymbol{\theta}_{0}, q_{0})$; (ii) consistency for the marginal posterior of the parametric component $\boldsymbol{\theta}_{0}$, after marginalizing out the nonparametric component $q$; or (iii) consistency for the marginal posterior of the nonparametric component $q_{0}$ after marginalizing out the parametric component $\boldsymbol{\theta}_{0}$. It is often the case that the nonparametric component is not of interest for inference, and is treated as a nuisance parameter. The remarkable results of \citet{ChengKosorok:2008} on the posterior profile distribution are of interest in that case.

\citet{Barron:1986}, \citet{GGR:1999b} and \citet{AGGR:2003} noted the importance of the Kullback-Leibler property in establishing consistency for the marginal posterior of the parametric component in semiparametric regression models. There have been many recent developments in Bayesian asymptotics for semiparametric models, a representative sample of which are summarized in \citet{CastilloRousseau:2015}, \citet{RSS:2014} and \citet{CSS:2016}.

\subsection{Rate of Posterior Contraction}

An implication of posterior consistency is the existence of a \emph{rate} of consistency, i.e. a sequence $\epsilon_{n}$ indexed by the sample size, corresponding to the size of a shrinking ball around the true density whose posterior probability tends to 1 as $n \rightarrow \infty$. We use the terms rate of \emph{contraction}, rate of \emph{convergence}, and rate of \emph{consistency} interchangeably. It is especially important in the study of Bayes factor asymptotics that one can establish the rate of convergence of the posterior distribution at the true value. The convergence rate is the size $\epsilon_{n}$ of the smallest ball, centered at the true value, such that the posterior probability of this ball tends to one.
\begin{definition}[\citet{GLV:2008}]
\label{def:convergerate}
Suppose $\mathbf{Y}^{(n)}$ are an i.i.d. random sample from $p_{0}$, $d$ is a distance on the set of densities. The convergence rate of the posterior distribution $\Pi_{n}(\cdot | \mathbf{Y}^{(n)})$ at the true density $p_{0}$ is at least $\epsilon_{n}$, if as $n \rightarrow \infty$, $\epsilon_{n}\rightarrow 0$, and for every sufficiently large constant $M$,
\[
\Pi_{n}(f: d(f,p_{0}) > M \epsilon_{n} | \mathbf{Y}^{(n)}) \rightarrow 0
\]
in $P_{0}^{n}$-probability.
\end{definition}
The sequence $\epsilon_{n} \rightarrow 0$ is often referred to as the targeted rate, where the target is the known optimal rate for the estimation problem at hand. Note that this definition utilizes the notion that complements of neighborhoods of the true density have posterior probability tending to zero.

We point out that, in fact, it is the \emph{higher-order asymptotic} properties of Bayesian methods which play a key role in Bayes factor consistency. That is to say, if we consider consistency to be a first-order property, then the study of rates of consistency can be viewed as a first step towards Bayesian nonparametric higher-order asymptotics. Such analysis will undoubtedly be very different from Bayesian parametric higher-order asymptotics, such as the posterior expansions reviewed in \citet{Ghosh:1994}.

The prior will completely determine the \emph{attainable} rate of posterior contraction, which can be understood intuitively by imagining how quickly the posterior will contract as the thickness of the prior tails is increased or decreased \citep{MartinWalker:2016}. However, the \emph{optimal} rate depends on the smoothness of the underlying true density, which is in general not known. In the regression setting, if the smoothness of the mean regression function is known, the prior can be suitably adjusted so that the attainable and optimal contraction rates coincide. When the smoothness is unknown, the prior should be more flexible so that it can adapt to the unknown optimal contraction rate. Of primary importance for Bayes factor asymptotics are the \emph{relative attainable} rates of posterior contraction in each model under comparison.

The core ideas regarding posterior convergence rates in the general (not necessarily parametric) setting were developed in \citet{GGV:2000} and \citet{ShenWasserman:2001}. We also mention \citet{GGV:2000, GhosalVDV:2007, Castillo:2014, GineNickl:2011} and \citet{HRS:2015}. Informative and didactic treatments are found in \citet[\S 7.3]{GineNickl:2015}, \citet[\S 2.5]{Ghosal:2010}, \citet{GLV:2008} and \citet{WLP:2007}. Recent extensions include misspecified models \citep{KleijnVDV:2006, Lian:2009, Shalizi:2009}, non-Euclidean sample spaces \citep{BhattacharyaDunson:2010}, and conditional distribution estimation \citep{PDT:2013}. Pseudo-posterior convergence rates have recently been considered in \citet{MHW:2013}.

\subsection{Connections with Other Literature}

Arguably the most closely-related strand of literature to Bayes factor asymptotics is that concerning Bayesian adaptation. Foundational work in this area is due to \citet{BelitserGhosal:2003} and \citet{Huang:2004}. We also mention \citet{GLV:2003}, \citet{Scricciolo:2006}, \citet{LemberVaart:2007} and \citet{VaartZanten:2009}. The convergence of the ratio of posteriors, under the probability law of the true model, will turn out to be determined by two factors which are well-studied in the Bayesian adaptation literature. First, the smoothness of the true density determines the minimax rate of convergence that any estimator can achieve. Second, the smoothness properties of the models under consideration will determine how closely their respective posteriors track the optimal minimax convergence rate. For recent progress on, respectively, adaptive minimax density estimation and Bayesian adaptation, see \citet{GoldenshlugerLepski:2014} and \citet{Scricciolo:2015}.

Parallel literatures concerning goodness-of-fit tests, prequential analysis, code length, and proper scoring rules contain many ideas relevant to Bayes factor asymptotics. There has been some cross-fertilization with the Bayes factor consistency literature. The interested reader is referred to \citet{TCG:2010}, \citet{Dawid:1992}, \citet{Grunwald:2007} and \citet{DawidMusio:2015}, respectively, for key ideas and references. Many of the asymptotic arguments presented in these parallel literatures are similar to those found in the statistics literature, though there is perhaps more emphasis in the former on information theory. \citet{ClarkeBarron:1990} and \citet{Barron:1998} are great resources for Bayesian asymptotics and the relevant information theory. Those authors derive the distribution of the relative entropy (Kullback-Leibler divergence), and in the process find bounds for the ratio of marginal likelihoods. They also give useful decomposition identities involving the Kullback-Leibler divergence and ratio of marginal likelihoods. \citet{Zhang:2006} further demonstrates how ideas from information theory can simplify existing results on posterior consistency.

\section{Comparing Two Parametric Models}
\label{sec:parpar}

The study of Bayes factor consistency in the parametric setting is more appropriately dealt with using different tools, rather than folding it into the nonparametric setting utilizing the Kullback-Leibler property of the priors. As noted by \citet{WDL:2004}, a finite-dimensional parametric model will not possess the Kullback-Leibler property unless $p_{0}$ is known to belong to the assumed parametric family. Moreover, if the parameter space $\Theta_{l}$ is nested in $\Theta_{k}$, then if $p_{0}^{n} \in \mathcal{M}_{l}$, the Kullback-Leibler property holds for both prior distributions prescribed by $\mathcal{M}_{k}$ and $\mathcal{M}_{l}$; c.f. \citet{RousseauChoi:2012}.

Consider the Bayes factor,$BF_{kl}=m(\mathbf{y} | M_{k})/m(\mathbf{y} | M_{l})$, for comparing any two models $\mathcal{M}_{k}$ and $\mathcal{M}_{l}$. From the {\em basic marginal likelihood identity} (BMI) of \citet{Chib:1995},
and momentarily suppressing subscripts, we have for each particular model that
\begin{equation}
m(\mathbf{y})=\frac{f(\mathbf{y}|\boldsymbol{\theta})\pi(\boldsymbol{\theta})}{\pi(\boldsymbol{\theta} | \mathbf{y})}.
\end{equation}
This equation is an identity because the l.h.s. does not depend on $\boldsymbol{\theta}$. Since this holds for any $\boldsymbol{\theta}$, if we consider any two particular sequences of $\boldsymbol{\theta}$ values indexed by the sample size $n$, the resulting sequences of marginal likelihoods will be the same. This simplification will allow us to focus on the maximum likelihood estimator for simplicity. Taking logs of both sides of this identity, we have
\begin{equation}
\log m(\mathbf{y})=\log f(\mathbf{y} | \boldsymbol{\theta})+\log \pi(\boldsymbol{\theta})-\log \pi(\boldsymbol{\theta} | \mathbf{y}),
\end{equation}
which yields
\begin{align}
\label{parBF}
BF_{kl}&=\frac{f_{k}(\mathbf{y}|\boldsymbol{\theta}_{k})\pi_{k}(\boldsymbol{\theta}_{k})}{\pi(\boldsymbol{\theta}_{k} | \mathbf{y})} \big/ \frac{f_{l}(\mathbf{y}|\boldsymbol{\theta}_{l})\pi_{l}(\boldsymbol{\theta}_{l})}{\pi(\boldsymbol{\theta}_{l} | \mathbf{y})} =\exp\{\log \frac{f_{k}(\mathbf{y}|\boldsymbol{\theta}_{k})}{f_{l}(\mathbf{y} | \boldsymbol{\theta}_{l})}+\log \frac{\pi_{k}(\boldsymbol{\theta}_{k})}{ \pi_{l}(\boldsymbol{\theta}_{l})}-\log \frac{\pi(\boldsymbol{\theta}_{k} | \mathbf{y})}{\pi(\boldsymbol{\theta}_{l} | \mathbf{y})}\},
\end{align}
or
\begin{equation}
\label{parlogBF}
\log BF_{kl}=\log \frac{f_{k}(\mathbf{y}|\boldsymbol{\theta}_{k})}{f_{l}(\mathbf{y} | \boldsymbol{\theta}_{l})}+\log \frac{\pi_{k}(\boldsymbol{\theta}_{k})}{ \pi_{l}(\boldsymbol{\theta}_{l})}-\log \frac{\pi(\boldsymbol{\theta}_{k} | \mathbf{y})}{\pi(\boldsymbol{\theta}_{l} | \mathbf{y})} .
\end{equation}
Note that \eqref{parlogBF} is also an identity. Thus, if we can find any two sequences $\tilde{\boldsymbol{\theta}}_{k,n}$ and $\tilde{\boldsymbol{\theta}}_{l,n}$ of consistent estimators of $\boldsymbol{\theta}_{k}$ and $\boldsymbol{\theta}_{l}$, respectively, such that $(i)$ $\log BF_{kl} \rightarrow_{p} -\infty$ when $M_{l}$ is the true model, and $(ii)$ $\log BF_{kl} \rightarrow_{p} \infty$ when $M_{k}$ is the true model, then this will establish Bayes factor consistency.

As discussed below, one may also be interested in what happens when the model dimension grows with the sample size. Thorough treatment of such scenarios would entail a major extension of the proposed framework, and the priors could no longer be viewed as independent of $n$.

\subsection{Nested Models}

We say that $M_{l}$ is nested in $M_{k}$ if $\Omega_{l} \subset \Omega_{k}$, and for $\boldsymbol{\theta} \in \Omega_{l}$, $f_{l}(\mathbf{y} | \boldsymbol{\theta})=f_{k}(\mathbf{y} | \boldsymbol{\theta})$. Mathematically this means that $\Omega_{l}$ is isomorphic to a subset of $\Omega_{k}$.

It is conceptually convenient to adopt a narrow interpretation of nested linear regression models, in which
\begin{equation*}
M_{k}: \; \mathbf{y} = \mathbf{X}_{k}\boldsymbol{\beta}_{k} + \boldsymbol{\varepsilon}, \;\; \boldsymbol{\varepsilon} \sim N_{n}(\mathbf{0}, \sigma_{k}^{2}I_{n}) , \;\;\;
M_{l}: \; \mathbf{y} = \mathbf{X}_{l}\boldsymbol{\beta}_{l} + \boldsymbol{\varepsilon}, \;\; \boldsymbol{\varepsilon} \sim N_{n}(\mathbf{0}, \sigma_{l}^{2}I_{n})
\end{equation*}
and $M_{l}$ is nested in model $M_{k}$ in the sense that columns of $\mathbf{X}_{l}$ are a proper subset of the columns of $\mathbf{X}_{k}$.

We examine \eqref{parlogBF} term-by-term, with $\tilde{\boldsymbol{\theta}}_{k,n}$ and $\tilde{\boldsymbol{\theta}}_{l,n}$ equal to the respective maximum likelihood estimators of $\boldsymbol{\theta}_{k}$ and $\boldsymbol{\theta}_{l}$.

Consider the first term evaluated at the maximum likelihood estimators under each model. Depending on which model is correct, this converges to a central or non-central chi-square random variable. It is thus bounded in probability, i.e. $O_{p}(1)$. A thorough treatment of the asymptotic distribution of the likelihood ratio statistic in the parametric setting, including nested, strictly non-nested, non-nested but overlapping, as well as misspecified models, is given by \citet{Vuong:1989}. In particular, this term is asymptotically distributed as either a weighted sum of chi-square random variables or standard normal, depending on whether the models are nested, overlapping or strictly non-nested. The Wilks phenomenon occurs when the asymptotic null distributions of test statistics are independent of nuisance parameters (or nuisance functions).  The Wilks phenomenon arises here in the special case that the models are nested. When the true parameter is on the boundary of the parameter space, some additional complications arise \citep{SelfLiang:1987}.

Next, consider the middle term. As $n \rightarrow \infty$, the prior effect on the Bayes factor is of order $O(1)$. The prior should not depend on the sample size unless the dimension of the parameter space changes (grows) with the sample size.

Evaluating the final term at the maximum likelihood estimators under each model, we see a potential problem in that if the maximum likelihood estimator is consistent in both models, with the same rate of consistency, then the sequence of ratios of posteriors will diverge, as each posterior density becomes concentrated around the true parameter value (and hence the numerator and denominator both diverge). We therefore require that the two models are asymptotically distinguishable in some sense. When both models are correctly specified, i.e. correctly specified and nested, some condition is needed to prevent the larger model from being chosen over the smaller model. The conventional argument achieves this by approximation of the Bayes factor using the Bayesian Information Criterion (BIC) suggested by \citet{Schwarz:1978}. This involves the Laplace approximation for the marginal likelihood, and the ratio of these marginal likelihood approximations in the Bayes factor yields a penalty term for model dimension. It should be noted, however, that it is possible for the Bayes factor to be consistent even when the BIC is not, as shown by \citet{BGM:2003} in the setting that the model dimension grows with the sample size.
Key references for the standard BIC-approximation-based argument for consistency of the Bayes factor include \citet{KassVaidyanathan:1992, GelfandDey:1994, OHaganForster:2004, Ando:2010} and  \citet{Dawid:2011}. We sketch a proof through the lens of our decomposition.

In \eqref{parlogBF}, we deal with the log ratio of posteriors by using a Laplace approximation for each posterior, which is naturally similar to using the Laplace approximation for the marginal likelihood. The Laplace approximation for the posterior given by \citet{Davison:1986} is
\begin{equation}
\label{LaplaceApprox}
\pi(\theta | \mathbf{y})=\frac{f(\mathbf{y} | \boldsymbol{\theta})\pi(\boldsymbol{\theta}) n^{-p/2} | I(\boldsymbol{\theta}^{*})|^{1/2}}{f(y | \boldsymbol{\theta}^{*})\pi(\boldsymbol{\theta}^{*})(2 \pi)^{p/2}}\{1+O_{p}(n^{-1})\},
\end{equation}
where $I(\boldsymbol{\theta}^{*})$ is minus the $p \times p$ matrix of second partial derivatives of $\log f(\mathbf{y} | \boldsymbol{\theta})+\log \pi(\boldsymbol{\theta})$ evaluated at the posterior mode $\boldsymbol{\theta}^{*}$.

The posterior mode can be replaced by the MLE since, asymptotically, the likelihood and its derivatives are negligible outside of a small neighborhood of the MLE. Plugging \eqref{LaplaceApprox} into \eqref{parlogBF} yields
\begin{align*}
\log BF_{kl}&=\log \frac{f_{k}(\mathbf{y}|\boldsymbol{\theta}_{k})}{f_{l}(\mathbf{y} | \boldsymbol{\theta}_{l})}+\log \frac{\pi_{k}(\boldsymbol{\theta}_{k})}{ \pi_{l}(\boldsymbol{\theta}_{l})}-\log \frac{\pi(\boldsymbol{\theta}_{k} | \mathbf{y})}{\pi(\boldsymbol{\theta}_{l} | \mathbf{y})}   \\
&= \log \frac{f_{k}(\mathbf{y}|\boldsymbol{\theta}_{k})}{f_{l}(\mathbf{y} | \boldsymbol{\theta}_{l})}+\log \frac{\pi_{k}(\boldsymbol{\theta}_{k})}{ \pi_{l}(\boldsymbol{\theta}_{l})} -\log\frac{f(\mathbf{y} | \boldsymbol{\theta}_{k})}{f(\mathbf{y} | \boldsymbol{\theta}_{k}^{*})} - \log \frac{\pi(\boldsymbol{\theta}_{k})}{\pi(\boldsymbol{\theta}_{k}^{*})} + \frac{(p_{k}-p_{l})}{2}\log n \\
& \;\;\; - \log \frac{| I(\boldsymbol{\theta}_{l}^{*})|^{1/2}}{| I(\boldsymbol{\theta}_{k}^{*})|^{1/2}} + \frac{(p_{k}-p_{l})}{2} \log (2\pi ) + \log\frac{f(\mathbf{y} | \boldsymbol{\theta}_{l})}{f(\mathbf{y} | \boldsymbol{\theta}_{l}^{*})} + \log \frac{\pi(\boldsymbol{\theta}_{l})}{\pi(\boldsymbol{\theta}_{l}^{*})}.
\end{align*}

Recall that $p_{l} < p_{k}$ since $M_{l}$ is nested in $M_{k}$. Now, simply take advantage of the identity and use the sequence of MLEs, $\hat{\boldsymbol{\theta}}_{k,n}=\boldsymbol{\theta}_{k,n}^{*}$. In this case the third term on the second line and the third term on the last line are both zero. Also the fourth terms on each line would be zero. Multiply what remains by minus two; the first term will be the usual $\chi_{p_{l}-p_{k}}^{2}$ statistic, and so is bounded in probability. The second term will be $O(1)$ (plug in the MLEs for both). The ratio of the determinants and the term involving $\pi$ are also $O(1)$. All that remains is $(p_{k}-p_{l}) \log n$. It is clear that for fixed-dimensional models, $-2 \log BF_{kl} \rightarrow -\infty$ in $P_{0}^{n}$-probability as $n \rightarrow \infty$.

\citet{MBR:1998} consider both nested and non-nested models using intrinsic Bayes factors, and establish consistency for nested models by arguing that the intrinsic Bayes factor approaches the usual Bayes factor in the limit. \citet{WangSun:2014} consider consistency of the Bayes factor for nested linear models when the model dimension grows with the sample size. \citet{MGC:2010,MGC:2015} study Bayes factor consistency in the same setting, with emphasis on objective Bayes factors, such as those using intrinsic priors. The intrinsic Bayes factors are shown to be asymptotically equivalent to the BIC, though not using the Laplace approximation. When the number of parameters is $O(n^{\upsilon})$ for $\upsilon <1$, the Bayes factor is consistent. The case $\upsilon=1$ is termed \emph{almost consistent} in the sense that consistency holds for all but a small set of alternative models. This small set of alternatives is expressed in terms of the pseudo-distance from the alternative model to the null model, as defined by \citet{MorenoGiron:2008}. \citet{ShangClayton:2011} consider models with growing $p$ in high-dimensional settings, and elucidate the connections with model selection consistency. \citet{JohnsonRossell:2012} considered Bayesian variable selection in high-dimensional linear models and found that the posterior probability of the true model tends to zero when the number of covariates is $O(n^{1/2})$, local priors (which assign prior density zero to null values) are used for regression coefficients, and the relative prior probabilities assigned to all models are strictly positive. This does not contradict Bayes factor consistency, which is based on pairwise comparisons.

\subsection{Non-Nested Models}

In the non-nested case, models may be either overlapping or strictly non-nested. Overlapping models are those for which neither model is nested in the other, but that there is at least one set of parameter values such that the conditional density of the observations is the same in the two models. In strictly non-nested models, there are no parameter values for which the conditional density of the observations is the same in two models. It is known that the BIC is generally not a consistent model selection criterion when selecting among non-nested models \citep{SinWhite:1996, HongPreston:2012}.

The standard treatment of non-nested linear regression models involves rival sets of predictor variables. The arguments for the non-nested linear regression model setting mirror those from the nested setting, with some important differences. The first term in \eqref{parlogBF} will converge to a different distribution; see \citet{Vuong:1989}. Provided that the two models are asymptotically distinguishable, only one will be correctly specified. The misspecified model posterior will not be consistent in the conventional sense, though there will still be a sort of pseudo-consistency, i.e. consistency for all pseudo-true parameter values, under certain conditions. See \S~\ref{sec:discussion} for key references. The log ratio of posterior densities will thus converge as desired, depending on which model is correctly specified. When both models are misspecified, some additional prior support conditions would be needed to ensure that the rate of posterior contraction is faster in the `best' model.

Few papers have considered this setting. \citet[p.1208]{CGMM:2009} state, ``As far as we know, a general consistency result for the Bayesian model selection procedure for non-nested models has not yet been established." Typically authors deal with the non-nested setting by specifying an encompassing model for which both models being compared are nested in the larger encompassing model. Let $\mathcal{M}_{\text{full}}$ be the full model including all predictors. This can be used as a base model to which $\mathcal{M}_{k}$ and $\mathcal{M}_{l}$ may be compared using
\[
BF_{kl}=\frac{BF_{k,\text{full}}}{BF_{l,\text{full}}},
\]
under a suitable choice of priors; see \citet{LPMCB:2008} and \citet{GuoSpeckman:2009}.  \citet{CGMM:2009} and \citet{GuoSpeckman:2009} establish consistency of Bayes factors for the special setting of normal linear models and a wide class of prior distributions, including intrinsic priors. The method of proof also utilizes the BIC approximation, though not the usual Laplace approximation (due to a property of intrinsic priors). \citet{WangMaruyama:2015} study consistency in linear models when the model dimension grows with the sample size, and argue that use of Zellner's $g$-prior is crucially important to establishing consistency. \citet{GMCM:2010} study objective Bayes factor consistency for non-nested linear models.

\section{Parametric vs. Nonparametric}
\label{sec:parnonpar}
\subsection{Nested Models}
A parametric model can be nested in a nonparametric model in several ways, of which the most common are: (i) the parametric model is a finite-dimensional restriction of the infinite-dimensional parametric model, (ii) in a regression setting, the set of predictors in the parametric model is a subset of the predictors in the nonparametric model. A non-nested example is the regression setting where the set of predictors is not the same in each model.

\subsubsection{The First Term}
The first term can be studied in at least two ways. The first is to consider a suitable probability inequality, perhaps derived from the results of \citet{WongShen:1995}, under which this term can be shown to be exponentially small. A second approach is to use the generalized likelihood ratio (GLR) results of \citet{FZZ:2001} where the Wilks phenomenon is seen to arise for suitable choices of the nonparametric density estimator. In particular, the local linear density estimator of \citet{Fan:1993}. Other important results with this theme include \citet{Portnoy:1988} and \citet{Murphy:1993}.

Generalized likelihood ratios were first studied by \citet{SeveriniWong:1992}. More recently, a detailed theory has been developed in \citet{FZZ:2001}, \citet{FanZhang:2004}, and \citet{FanJiang:2005, FanJiang:2007}. Consider a vector of functions, $\mathbf{p}$ and $\boldsymbol{\zeta}$, which are the parameters of a semiparametric or nonparametric model. Given $\boldsymbol{\zeta}$, a nonparametric estimator of $\mathbf{p}$ is given by $\hat{\mathbf{p}}_{\boldsymbol{\zeta}}$. The parameters $\boldsymbol{\zeta}$ are regarded as nuisance parameters and are estimated using the profile likelihood, i.e. by finding $\boldsymbol{\zeta}$ to maximize $\ell(\hat{\mathbf{p}}_{\boldsymbol{\zeta}}, \boldsymbol{\zeta})$ with respect to $\boldsymbol{\zeta}$. This yields the maximum profile likelihood $\ell(\hat{\mathbf{p}}_{\hat{\boldsymbol{\zeta}}}, \hat{\boldsymbol{\zeta}})$. \citet{FanJiang:2007} emphasize that this is not a maximum likelihood, since $\hat{\mathbf{p}}_{\boldsymbol{\zeta}}$ is not an MLE. Suppose we are testing a parametric null hypothesis against a nonparametric alternative, i.e.
\[
H_{0}: p = p_{\theta}, \; \theta \in \Theta .
\]
Denote the MLE under the null model as $(\hat{\theta}_{0}, \hat{\boldsymbol{\zeta}}_{0})$, which maximizes the log-likelihood $\ell(\mathbf{p}_{\theta}, \boldsymbol{\zeta})$. Then $\ell(\mathbf{p}_{\hat{\theta}_{0}}, \hat{\boldsymbol{\zeta}}_{0})$ is the MLE under the null. The GLR statistic is
\[
GLR_{n}=\ell(\hat{\mathbf{p}}_{\hat{\boldsymbol{\zeta}}}, \hat{\boldsymbol{\zeta}})-\ell(\mathbf{p}_{\hat{\theta}_{0}}, \hat{\boldsymbol{\zeta}}_{0}).
\]
What is amazing about the GLR theory is that it is not necessary to use a genuine likelihood; quasilikelihoods may also be employed. In broad generality, the asymptotic null distribution of the GLR statistic is approximately $\chi^{2}$.

A particularly attractive feature of the GLR framework in our study is that it allows for flexibility in the choice of estimator used for the nonparametric model. The method of estimation and the smoothing parameters used will in general affect $\mu_{n}$ and $r$, so that the relevant asymptotic distribution would need to take the estimation method into account. However, the basic result of convergence in distribution is broadly applicable, and this alone suffices for our purposes. In particular, we can choose any convenient nonparametric estimator to achieve the desired rate of posterior contraction needed for Bayes factor consistency, while remaining confident that any reasonable choice will ensure convergence in distribution of the GLR statistic.

\subsubsection{The Third Term}

We need some additional assumptions for this problem to be well-defined. Here we mention only the assumptions which are most important for intuition, but warn the reader that other technical conditions can be found in the cited papers that are essential for proving Bayes factor consistency. In particular, suppose that \citep{GLV:2008, MRM:2009} $\mathcal{M}_{k}$ prescribes a finite-dimensional parametric set of densities $\mathcal{F}_{k}$.
\begin{assumption}
\label{PvNP_A1}
Let $\epsilon_{n}$ be a sequence of numbers such that $\epsilon_{n} \rightarrow 0$ as $n \rightarrow \infty$. The nonparametric posterior $\Pi_{l}(\cdot | \mathbf{Y}^{(n)})$ is strongly consistent at $p_{0}$ with rate $\epsilon_{n}$ in the sense that
\[
\Pi_{l}(f : d(f,p_{0})>\epsilon_{n} |\mathbf{Y}^{(n)}) \rightarrow 0
\]
almost surely, in $P_{0}^{\infty}$-probability, where $d$ is some distance on the space of densities such as Hellinger or $L_{1}$.
\end{assumption}
This assumption is satisfied by most commonly-used nonparametric priors. If the posterior is not consistent, the model comparison problem would not be well-defined. In existing methods of proof, it is common to also make an assumption which ensures that the marginal likelihood under the parametric model is bounded from below.
\begin{assumption}
Let $f_{\theta}$ denote a generic element of the parametric family prescribed by $\mathcal{M}_{k}$. For any $\theta \in \Theta \subset \mathbb{R}^{s}$, with $s$ finite,
\[
\Pi_{k}(\theta^{\prime} : d_{KL}(f_{\theta}, f_{\theta^{\prime}})<cn^{-1}, \; V(f_{\theta},f_{\theta^{\prime}})<cn^{-1})>Cn^{-s/2},
\]
with $V(p,q)=\int p \log(p/q)^{2}d\mu$, where $p,q$ are both absolutely continuous with respect to $\mu$, $c$ and $C$ are positive constants.
\end{assumption}
Most models which for which a formal Laplace expansion is valid will satisfy this assumption. If both models are correctly specified, then to ensure that the simpler, parametric model is chosen, the rate of posterior contraction in the parametric model must be sufficiently faster than the corresponding rate of posterior contraction in the nonparametric model.
\begin{assumption}
Let $A_{\epsilon_{n}}(\theta)=\{ f: d(f,f_{\theta})<C \epsilon_{n}\}$, where $\epsilon_{n}$ is the rate of consistency of the nonparametric posterior from Assumption~\ref{PvNP_A1}. Then
\[
\sup_{\theta}\Pi_{l}(A_{\epsilon_{n}}(\theta))=o(n^{-d/2}) .
\]
\end{assumption}
This controls the amount of mass that the nonparametric prior assigns to such neighborhoods of the parametric family of densities, for all values of $\theta \in \Theta$. Given this assumption and the rate of consistency for the nonparametric posterior, then the Bayes factor is consistent.

\subsection{Non-Nested Models}

When the parametric model is not nested in the alternative, the GLR framework does not help with the first term. In that case it may be easiest to utilize relevant probability inequalities to bound the log-likelihood ratio statistic so that it is exponentially small with probability exponentially close to 1. The third term can be dealt with similarly to the nested case, as discussed in detail in \citet{GLV:2008}[\S 4].

\subsection{Connections to Literature}
Bayes factor consistency in finite-dimensional versus infinite-dimensional model comparison problems is addressed for i.i.d. data in \citet{GLV:2008} and \citet{MRM:2009}. The assumptions in the former paper are similar in spirit to those in the latter, which are given above. The results discussed below regarding the independent but not identically distributed setting for comparing a linear and partially linear model \citep{ChoiRousseau:2015} are also relevant here.

Much of the existing literature in this area considers such model comparisons as a goodness of fit testing problem. An excellent review is found in \citet{TCG:2010}. \citet{DassLee:2004} considered testing a point null versus a nonparametric alternative. A common approach is to embed a parametric model into a nonparametric alternative using, for example, a mixture of Dirichlet processes \citep{CarotaParmigiani:1996}, a mixture of Gaussian processes \citep{VerdinelliWasserman:1998}, or a mixture of Polya tree processes \citep{BergerGuglielmi:2001}. These methods are connected to the approach of \citet{FRR:1996}, in which a Dirichlet process prior was used for the alternative. Recent efforts in this area are discussed in \citet{AlmeidaMouchart:2014}.

We mention that one well-studied example in this setting is the comparison of a nonparametric logspline model with a parametric alternative. The asymptotic properties of the nonparametric logspline model have been studied in \citet{Stone:1990, Stone:1991}. \citet{JWC:2010} study the consistency property of information-based model selection criteria for this model comparison problem. \citet{GLV:2008} use the logspline model as an example to illustrate that the prior support conditions are satisfied, in order to prove Bayes factor consistency in this setting.

\section{Parametric vs. Semiparametric}
\label{sec:parsemi}

It is increasingly common in regression modeling to compare parametric linear models against semiparametric models. A typical proof of Bayes factor consistency in this setting is found in \citet{KunduDunson:2014}. The proof involves bounding the marginal likelihood (upper and lower bounds) in each model and then considering what happens as $n \rightarrow \infty$. The arguments closely follow those of \citet{GuoSpeckman:2009}, and can also accommodate increasing-dimensional parameter spaces, subject to a condition on the rate of growth in terms of the available sample size.

As in the parametric versus nonparametric comparison problem, existing results regarding posterior contraction rates in specific problems can be adapted to the study of the log ratio of posteriors in our framework. It is the log ratio of likelihoods which requires some additional care.

The log likelihood ratio can often be handled by some variant of semiparametric likelihood ratio tests. The results of \citet{MurphyVDV:1997} are particularly useful when the finite dimensional parametric component is of interest, and the infinite-dimensional nonparametric component is treated as a nuisance parameter. In some cases, the restricted likelihood ratio test may be used, as described in \citet{RWC:2003}[Ch. 6]. That particular formulation relies on modeling the mean regression function as a penalized spline having a mixed effects representation, and is potentially useful when comparing a parametric linear mean regression model against an encompassing semiparametric model. This approach can also handle some types of dependent data; see the pseudo-likelihood ratio test proposed by \citet{SLCR:2014}.

In the semiparametric regression setting when one is interested in inference for the nonparametric component, a promising approach is described by \citet{LiLiang:2008}. Those authors extend the generalized likelihood ratio test described in \S~\ref{sec:parnonpar} to semiparametric models, which they use for inference about the nonparametric component. Remarkably, the Wilks phenomenon is again seen to arise for the generalized varying-coefficient partially linear model, which includes partially linear models and varying coefficient models as special cases.

Dealing with the log ratio of posteriors is of similar difficulty. When a parametric null model is being compared to a partially linear model, such that the null model is nested in the alternative, and if the true density is contained in the null model, it is clear that both models will possess the Kullback-Leibler property. This situation bears resemblance to the parametric versus nonparametric model comparison problem already discussed. For i.i.d. observations, conditions for Bayes factor consistency are studied in \citet{Rousseau:2008}, \citet{GLV:2008} and \citet{MRM:2009}. The case of independent but not identically distributed observations is studied in \citet{ChoiRousseau:2015}.
\begin{example}
\label{partlinear_compare}
\begin{equation}
\mathcal{M}_{1}: \; Z_{i}= W_{ij} \boldsymbol{\beta} + \sigma \varepsilon_{i}, \;\;\; \mathcal{M}_{2}: \; Z_{i}= W_{ij} \boldsymbol{\beta} + r(X_{i})+ \sigma \varepsilon_{i},
\end{equation}
where $r(\cdot)$ is an unknown function in an infinite-dimensional parameter space, and $W_{ij}\boldsymbol{\beta}$ is a linear component.
\end{example}
A particular version of this comparison problem is studied by \citet{ChoiRousseau:2015}, where $\{W_{ij}\} \in [-1,1]$, $i=1, \ldots, n$, $j=1, \ldots, p$, with $p$ finite, and $\{X_{i}\} \in [0,1]$, $i=1, \ldots, n$. Recall that a semiparametric prior is $\Pi = \Pi_{\boldsymbol{\beta}} \times \Pi_{r}$. \citet{ChoiRousseau:2015} assume a Gaussian process prior for $\Pi_{r}$ with a reproducing kernel Hilbert space (RKHS) $\mathbb{H}=L^{2}([0,1])$, and concentration function $\phi_{r}$ given by $\phi_{r}(\epsilon) = \inf_{h \in \mathbb{H}: \|h-r_{0} \|_{\infty} < \epsilon} \|h \|_{\mathbb{H}}^{2}-\log \Pi_{r} [ \|r \|_{\infty} < \epsilon ]$, where $\| \cdot \|_{\infty}$ is the supremum norm, $\| r \|_{\infty} = \sup \{ |r(x) | : x \in [0,1]\}$. Concentration functions and RKHS are discussed in \citet{VaartZanten:2008a,VaartZanten:2008b}. They show that the support of $\Pi_{r}$ is the closure of $\mathbb{H}$ and, therefore, for any $r_{0}$ contained in the support of $\Pi_{r}$, there exists $\epsilon_{n} \downarrow 0$, and $c>0$ such that $\phi_{r_{0}}(\epsilon_{n}) \leq n \epsilon_{n}^{2}$, and the following small ball probability holds
\[
\Pi_{r}[\|r-r_{0}\|_{\infty} \leq \epsilon_{n}] \geq e^{-cn \epsilon_{n}^{2}}.
\]
This property ultimately controls the rate of posterior contraction under the partially linear model, and is essential to the proof of Bayes factor consistency in \citet{ChoiRousseau:2015}.

\citet{ChoiRousseau:2015} and \citet{KunduDunson:2014} give regularity conditions for Bayes factor consistency in the semiparametric setting, both for convergence in probability and convergence almost surely. \citet{Lenk:2003} discusses testing a parametric model against a semiparametric alternative by embedding an exponential family in a semiparametric model, and develops an MCMC procedure for evaluating the adequacy of the parametric model. \citet{CLR:2009} study the asymptotic properties of the Bayes factor when comparing a parametric null model to a semiparametric alternative. They impose strong conditions on the structure of the model, namely Gaussian errors with known variance, and an orthogonal design matrix, and a particular variance structure for the coefficients in a trigonometric series expansion representation of the nonparametric component of the mean function in the regression model. An interesting example of this model comparison setting is found in \citet{MacEachernGuha:2011}, which illustrates a seemingly paradoxical result that under certain prior support assumptions, the posterior under a semiparametric model can be more concentrated than under a finite-dimensional restriction of that model.

A different strand of semiparametric Bayesian literature exists for models defined by moment conditions. Empirical likelihood methods have become popular for frequentist estimation and inference in such models. Empirical likelihoods have many properties of parametric likelihoods, but until recently there was no formal probabilistic interpretation, which meant there was no Bayesian justification for their use. \citet{Schennach:2005, Schennach:2007} introduced the Bayesian exponentially-tilted empirical likelihood (ETEL) to study the marginal posterior of a finite-dimensional interest parameter in the presence of an infinite dimensional nuisance parameter. This analysis is extended by \citet{CSS:2016} to the problem of selecting valid and relevant moments among different moment condition models. Using a semiparametric approach, it is shown that the ETEL satisfies a Bernstein-von Mises theorem in misspecified moment models. Moreover, the marginal likelihood-based moment selection procedure, based on the approach of \citep{Chib:1995}, is proven to be consistent in the sense of selecting only correctly specified and relevant moment conditions. This significant new work opens up the Bayes
factor computation and analysis for a large array of Bayesian semiparametric models that are specified only through the moments of the underlying unknown distribution. It also connects to the analysis of similar problems from a frequentist perspective such as the proportional likelihood ratio model \citep{LuoTsai:2012}, which can be seen as a generalization of exponential tilt regression models, and its adaptation to mean regression modeling using empirical likelihood \citep{HuangRathouz:2012}. These frequentist developments have spawned new results concerning likelihood ratio tests. For example, in the context of high-dimensional semiparametric generalized linear models, the reader is referred to \citet{NZL:2015}.

\section{Directions for Future Research}
\label{sec:discussion}

We briefly discuss some potential avenues for future research which coincide with recent advances in the literature.

\textbf{Other Model Comparison Problems.} Some model comparison problems are less well-studied than the three model-comparison frameworks already discussed: (i) two nonparametric models; (ii) two semiparametric models; and (iii) a nonparametric and a semiparametric model. In principle, the concepts we have discussed in other cases continue to guide our intuition and direct our approach to the problem. However, we are unaware of existing theoretical results which are general enough to be broadly applicable in these settings.  An ongoing area of research is that of testing nonparametric hypotheses when they are not well-separated; see \citet{Salomond:2015}.

\textbf{Rate of Convergence of Bayes Factors.} \citet{MRM:2009} actually obtain an upper bound on the convergence rate of the Bayes factor when the parametric model is correct and the comparison is with an encompassing nonparametric model. Such results would be of great use in applied settings. In particular, it would be helpful to be able to state more precisely the notion that under certain prior support conditions, the Bayes factor cannot distinguish between a correct and incorrect model, or may even prefer an incorrect model, despite having very large sample sizes.

\textbf{Empirical Bayes.}  Replacing some of the hyperparameters in the prior distribution by data dependent quantities is part of an empirical Bayesian analysis. The asymptotic properties of empirical Bayes procedures have recently been studied in nonparametric models, where the maximum marginal posterior likelihood estimator is used for the hyperparameters; see \citet{PRS:2014} and \citet{RousseauSzabo:2015}. Consistency of the posterior in empirical Bayesian analysis is complicated by the need to also consider the sequence of estimators for the hyperparameter. Posterior concentration rates with empirical priors are studied in \citet{MartinWalker:2016}.

\textbf{Summary Statistics.}  An exciting new approach to Bayes factor asymptotics by \citet{MPRR:2014} shows how Bayes factor consistency can hinge on whether or not the expectations of suitably chosen summary statistics are asymptotically different under the models being compared. It will be interesting to learn if those conditions, which are necessary and sufficient, can be related in a meaningful way to conditions on the log ratio of posterior densities.

\textbf{Non-i.i.d. Observations and Conditional Densities.}  A natural extension to existing results concerns non-i.i.d. observations. Consistency and asymptotic normality for posterior densities in such settings are active research areas. See \citet{GhosalVDV:2007}. A related extension is the estimation of conditional densities. For example, in a semiparametric regression model, one might assume a linear model in the predictors and a nonparametric error distribution. One may not wish to assume that the predictors and errors are independent, and then could consider inference for the nonparametric component, i.e. the conditional error density, given the predictors. See e.g. \citet{Pelenis:2014}. \citet{PDT:2013} investigate the possibility of recovering the full conditional distribution of the response given the covariates.

\textbf{Misspecified Models.}  Another important question is the large-sample behavior of the posterior when the model is misspecified. It is shown in \citep[p.39-40]{vanOmmen:2015} that the posterior is inconsistent in simple nested model settings with misspecification (the truth is not in the set of models being considered) despite all the relevant consistency theorems still holding, e.g. \citet{KleijnVDV:2006}, \citet{DeBlasiWalker:2013}, and \citet{RSM:2015}. An enlightening discussion of some of these issues, with more references, is given by \citet{VGMRW:2015}. Other key references include \citet{BunkeMilhaud:1998}, \citet{GrunwaldLangford:2007}, \citet{Shalizi:2009}, \citet{LeeMacEachern:2011}, \citet{Muller:2013}, and \citet{CSS:2016}. Section 4.3.2 of a recent Ph.D. thesis \citep{vanOmmen:2015} contains a current discussion of Bayesian consistency under misspecification. \citet{OSS:2015} consider the lack of robustness of Bayesian procedures to misspecification, which they call `brittleness'. \citet{KleijnVDV:2012} study the Bernstein-von Mises theorem under misspecification as do \citet{CSS:2016}.

\end{document}